# Dyck path factorization

Gennady Eremin
ergenns@gmail.com

November 26, 2021

**Abstract.** Dyck paths are among the most heavily studied Catalan families. We work with peaks and valleys to uniquely decompose Dyck paths into the simplest objects − *prime fragments* with a single peak. Each Dyck path is uniquely characterized by a set of peaks or a set of valleys. The appendix contains a python program with which the reader can factorize Dyck paths online.

*Keywords*: factorization, decomposition, Dyck path, Dyck word, Catalan number, prime fragment, height, Cantor's pairing function.

## 1 Introduction

Dyck paths are among the most heavily studied Catalan families [1]. There are many structures equivalent to Dyck paths, all of which are counted by the Catalan numbers (see A000108 in [11]). In the literature, there is a natural decomposition of the Dyck path into paired nested up steps and down steps (see [2], p. 169). A similar decomposition takes place in the case of well-matched opening and closing parentheses. However, as a result, the number of decomposition objects is commensurate with the length of the Dyck path (or the length of the well-parenthesed word).

Another approach, factorization is often understood as the decomposition of a given Dyck path into fragments that are themselves Dyck paths (see [3], p. 3). Concatenation of these fragments gives the original Dyck path. As it often happens, the best results are obtained by combining different methods.

Quite a lot of papers explore peaks and valleys in Dyck paths for various purposes (e.g., [4 − 8]). We will work with peaks and valleys to uniquely decompose Dyck paths into the simplest objects − *prime fragments* with a single peak. In this paper we use some results obtained in [9, 10].

## 2 Dyck path terminology

We will adhere to the terminology [2, 9]. A Dyck path of semilength *n* is a diagonal lattice path in the first quadrant with up steps $u = \langle 1, 1 \rangle$, *rises*, and down steps $d = \langle 1, -1 \rangle$, *falls*, that starts at the origin (0, 0), ends at (2n, 0), and never passes below the *x*-axis. The Dyck path of semilength *n* we will call an *n*-Dyck path.

If in a given Dyck path we encode each rise with the letter *u* and each fall with the letter *d*, then we get the corresponding *Dyck word* over the alphabet {*u*, *d*} in which there are equal numbers of occurrences of the letters *u* and *d*, with at least as many occurrences of the letter *u* in any initial subword as the letter *d*.

A point of a Dyck path with ordinate *y* is said to be at *level y*. By a *return step* we mean a down step with an endpoint at the level 0 (the *x*-axis). In a Dyck path, a *peak* (lo-cal maximum) is the inner point of a subword of the form *ud*, a *valley* (local minimum) is the inner point of a subword of the form *du*. A maximum string of *u*'s (*d*'s) is an *ascent* (*descent*) of a Dyck path.

Below, Figure 1 shows the 6-Dyck path *ud-uuudd-uuddd* that we borrowed from [2] (the dash marks inside the Dyck word will be explained a little later). As you know, all *n*-Dyck paths are bounded by the *n-Dyck triangle* with three vertices $(0, 0)$, $(n, n)$ and $(2n, 0)$ (a pyramid of height *n* that corresponds to word $u^n d^n$).

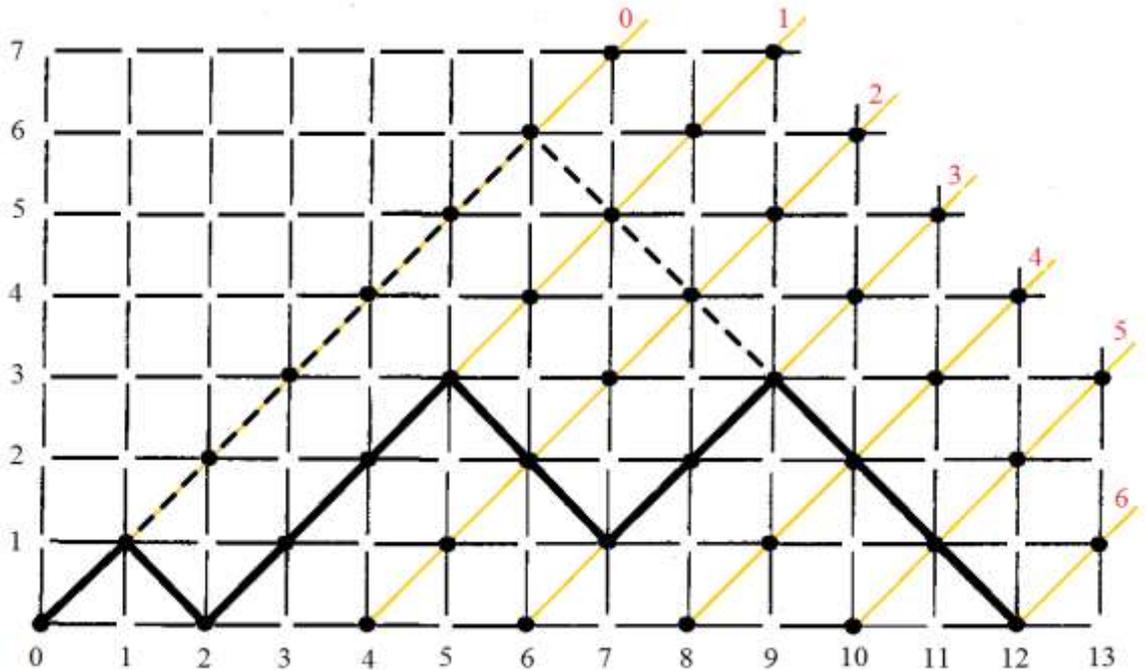

Figure 1. The 6-Dyck path *ud-uuudd-uuddd*.

In the figure, we have marked with dots non-negative lattice nodes achievable for Dyck paths. As we can see, not all nodes are reachable. Anything above the main diagonal $y = x$ is not available. Below the main diagonal, only nodes with an even sum of coordinates are reachable. With a different color, we showed and numbered the ascending diagonals on which the up steps are placed. In order not to obscure the drawing, we did not show descending diagonals (for down steps).

**Definition 1.** *A prime fragment of a Dyck path is a peak together with an adjacent ascent (left) and descent (right).*

In accordance with Definition 1, the Dyck path with a single peak, a *Dyck triangle* [9], is called *prime*. In this case, the extreme vertices of a single prime fragment lie on the *x*-axis.

**Definition 2.** *Factorization of the Dyck path is its decomposition into prime fragments.*

In the given 6-Dyck path, three prime fragments are separated by a dash: (i) *ud* with peak $p_1 = (1, 1)$ inside, (ii) *uuudd* with peak $p_2 = (5, 3)$ inside, and (iii) *uuddd* with peak $p_3 = (9, 3)$ inside. Additionally, we note two valleys: $v_1 = (2, 0)$ and $v_2 = (7, 1)$. Conditionally, we can consider the extreme points of this path as valleys: $v_0 = (0, 0)$ and $v_3 = (12, 0)$.

It is not difficult to see, peaks and valleys strictly alternate in each $n$-Dyck path. These points can be ordered in ascending order $x$-coordinate from $(0, 0)$ to $(2n, 0)$. Obviously, a set of peaks and valleys is uniquely formed for the Dyck path. The converse is also true, for a given set of peaks and valleys, it is easy to reconstruct the Dyck path. Thus we have a bijection.

Next, we will show that in order to obtain a bijection, it is sufficient to select either only peaks or only valleys.

## 3 Peaks and valleys of the Dyck path

Let's take a closer look at the situation with neighboring peaks and valleys. At each peak of the Dyck path, the adjacent ascent and descent form a 90-degree angle. If we continue both sides of the angle to the $x$-axis, then we get an isosceles triangle in which the height divides the base in half and each half equal to the height.

**3.1. Adjacent peak, valley and peak.** Figure 2 shows a fragment of an $n$-Dyck path with adjacent peaks: left $p_l = (x_l, y_l)$ and right $p_r = (x_r, y_r)$, $y_l < y_r$. The valley $v = (x_v, y_v)$ is located between two peaks. Obviously, $x_l = y_l$ if $l = 1$. If the right peak $p_r$ is the last one, then $x_r + y_r = 2n$.

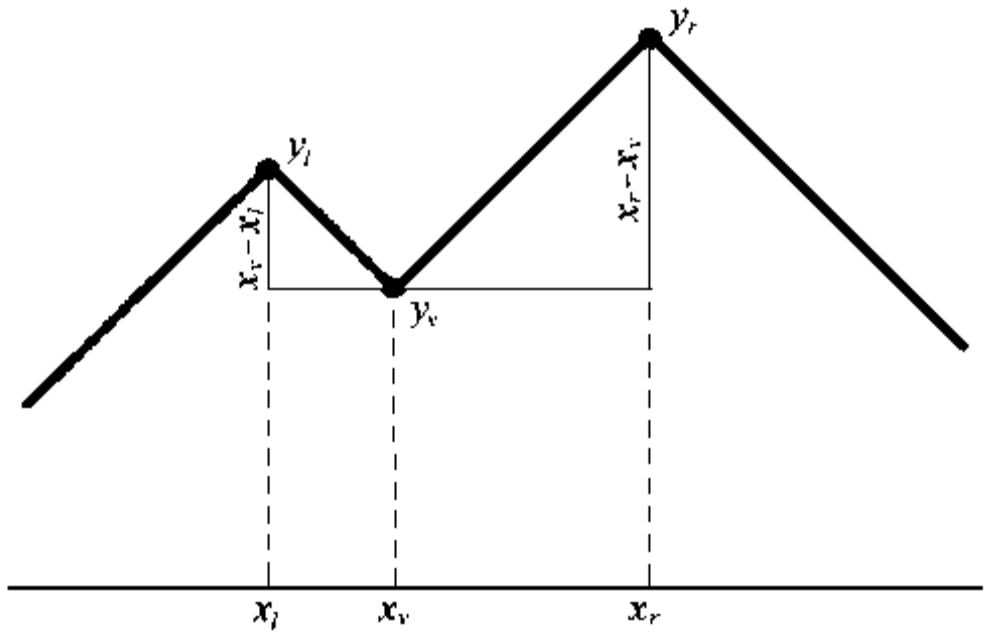

Figure 2. The valley between two peaks.

It is quite simple to calculate the coordinates of the valleys $v$:

$$x_v - x_l = y_l - y_v \text{ or } x_v + y_v = x_l + y_l,$$
$$x_r - x_v = y_r - y_v \text{ or } x_v - y_v = x_r - y_r.$$

If we denote $a = x_l + y_l$ and $b = x_r - y_r$, then we get

(1) $$x_v = (a + b)/2 \text{ and } y_v = (a - b)/2.$$

Obviously, always $a \leq b$. In the case $a = b$, the valley lies on the x-axis, i.e. $y_v = 0$.

Note that formula (1) is valid regardless of the height of neighboring peaks, both $y_l \leq y_r$ and $y_l > y_r$ are possible. Thus, the set of peaks of the Dyck path is enough to restore all the valleys and hence the entire Dyck path. For example, the 6-Dyck path in Figure 1 is represented by a set of three vertices $\{(1, 1), (5, 3), (9, 3)\}$. We formulate the result as a theorem.

**Theorem 1**. *The Dyck path is completely determined by the set of its peaks.*

**3.2. Adjacent valley, peak and valley.** The following Figure 3 shows a fragment of an $n$-Dyck path with adjacent valleys: left $v_l = (x_l, y_l)$ and right $v_r = (x_r, y_r)$, $y_l > y_r$. The peak $p = (x_p, y_p)$ is located between two valleys. Now we need to find the coordinates of the peak $p$. Here we can note two special cases: (i) if $x_l = y_l = 0$, then $x_p = y_p$; (ii) in the case of $x_r = 2n$ (respectively $y_r = 0$) we get $x_p + y_p = 2n$.

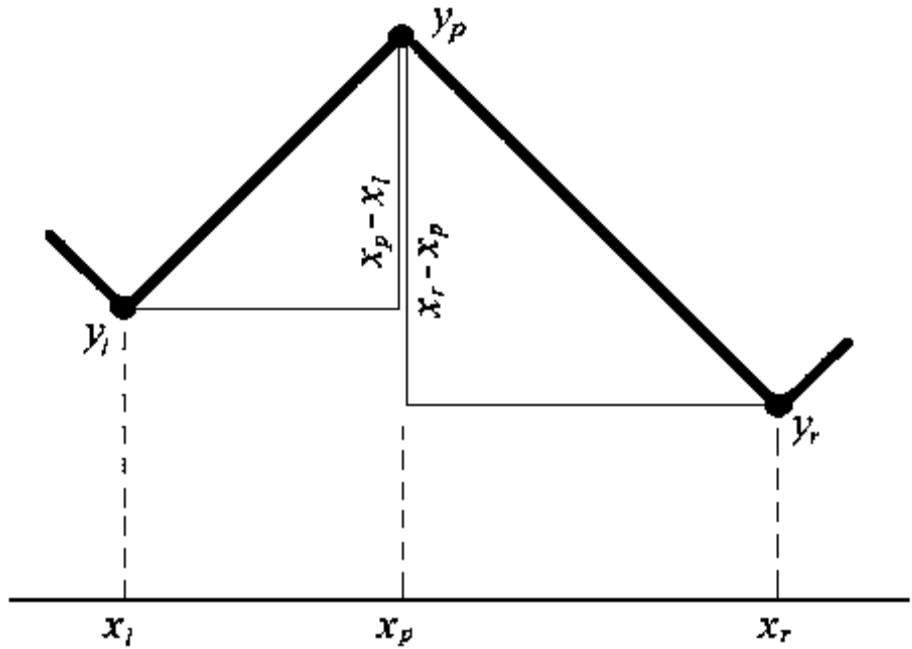

Figure 3. The peak between two valleys.

The following equalities follow directly from Figure 3:

$$x_p - x_l = y_p - y_l \quad \text{or} \quad x_p - y_p = x_l - y_l,$$
$$x_r - x_p = y_p - y_r \quad \text{or} \quad x_p + y_p = x_r + y_r.$$

Again we denote $s = x_l - y_l$ and $t = x_r + y_r$, as a result we get

(2) $\qquad\qquad x_p = (s + t)/2 \quad \text{and} \quad y_p = (t - s)/2.$

Thus, the set of valleys of the Dyck path is enough to restore all the peaks and hence the entire Dyck path. For example, the 6-Dyck path in Figure 1 is represented by a set of three valleys $\{(2, 0), (7, 1), (12, 0)\}$. Let's formulate another obvious theorem.

**Theorem 2**. *The Dyck path is completely determined by the set of its valleys.*

Dyck paths are often drawn in various coordinate grids. In Figure 1, the most common coordinate axes are used. However, in this case, only a quarter of the nodes of the first quadrant are available for Dyck paths. It can be said that unreachable grid nodes are ballast; the presence of inaccessible nodes complicates analysis and calculations.

## 4 Modified Dyck triangle

In [10], six different projections of the 4D Dyck triangle are considered. Figure 4 below shows the most compact version (see [10], Figure 6).

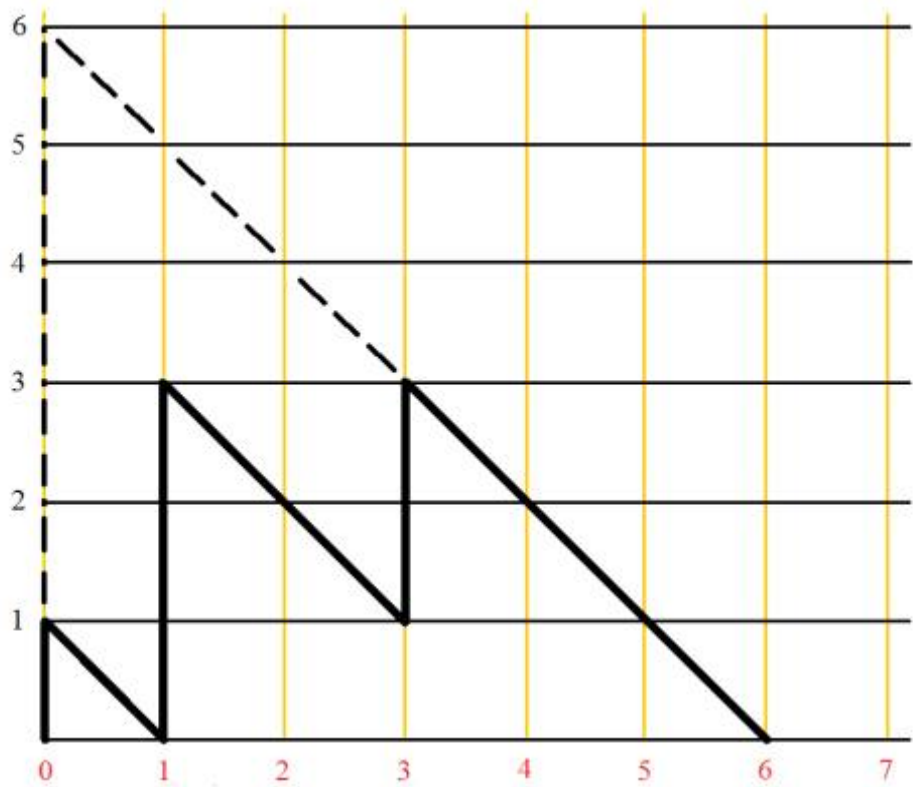

Figure 4. The modified 6-Dyck path.

The difference of this projection is that all points of the first quadrant are achievable for Dyck paths. In this case, we have condensed the *x*-axis using ascending diagonals from Figure 1 (shown in yellow). Each reachable node of the previous coordinate grid received a new *x*-coordinate $x' = (x - y)/2$. The y-coordinate has not changed.

All modified *n*-Dyck paths are bounded by a triangle with vertices $(0, 0)$, $(0, n)$ and $(n, 0)$. In a given 6-Dyck path the set of peaks is $\{(0, 1), (1, 3), (3, 3)\}$ and the set of valleys is $\{(1, 0), (3, 1), (6, 0)\}$.

## 5. Conclusion

The considered sets of peaks and valleys can be useful for numbering of Dyck paths (the topic of the next paper). To do this, we additionally need Cantor's pairing function, with which we can condense the resulting sets of integers [12, 13].

The appendix contains a python program with which the reader can factorize Dyck paths online (using Python IDE [14]). In the first part of the program, the received Dyck word is checked. The second part of the program builds a set of peaks for the modified word.

Gzhel State University, Moscow, 140155, Russia
*Email address*: http://www.art-gzhel.ru/

**Appendix.** The Python program calculates a set of peaks for the modified Dyck word.

```
word = ''; n2 = 0
while True:     # checking the Dyck word
    word = input('Please enter a Dyck word (alphabet: u, d)')
    n2 = len(word);
    if n2 == 0: continue
    y = 0;
    for x in range(n2):
        symb = word[x]; print(symb, end=' ')
        if symb == 'u': y += 1; continue
        if symb == 'd':
            y -= 1;
            if y >= 0: continue
        y = -1; break
    if y == 0: break
    print('- invalid word');
print('- Ok! Word length is', n2)
print('Peaks in the modified Dyck path: ', end='')

y = 0;    # calculating the peaks for the modified word
for x in range(n2):
    if word[x] == 'u':
        prior = 'u'; y += 1; continue
    if prior == 'u':
        print((x-y)//2, ',', y, '; ', end='')
    prior = 'd'; y -= 1
```